\documentclass{article}   	
\usepackage{geometry}               
\usepackage{graphicx}
\usepackage[usenames,dvipsnames]{xcolor}
\usepackage{amssymb,amsmath,amsthm,amsfonts}
\usepackage{bm}                     
\usepackage[english]{babel}
\usepackage[utf8]{inputenc}
\usepackage{wrapfig}
\usepackage[colorlinks]{hyperref}
\hypersetup{linkcolor=blue,citecolor=blue,filecolor=black,urlcolor=blue}

\usepackage[sc]{mathpazo}

\newtheorem{proposition}{Proposition}[section]
\newtheorem{theorem}[proposition]{Theorem}
\newtheorem{corollary}[proposition]{Corollary}
\newtheorem{lemma}[proposition]{Lemma}

\theoremstyle{definition}

\theoremstyle{remark}
\newtheorem{remark}[proposition]{Remark}
\numberwithin{equation}{section}
\newcommand{\abs}[1]{\left|#1\right|}
\newcommand{\norm}[1]{\left\lVert#1\right\rVert}
\newcommand{\epsen}[2]{\mathcal{E}_{\varepsilon}\left(#1, #2\right)}
\newcommand{\ep}{\varepsilon}
\newcommand{\eps}{\varepsilon}

\newcommand{\cf}{C_{f}}
\newcommand{\cl}{C_{L}}
\newcommand{\qc}{q_c}
\newcommand{\cpr}{C_{q}}
\newcommand{\leb}[1]{L^{#1}(\Omega)}
\newcommand{\rn}{\mathbb{R}^{N}}
\newcommand{\ml}[2]{#1L\left(-\frac{#2}{#1}\right)}
\newcommand{\nrg}[2]{\int_{\Omega} \ml{#1}{#2} + F(x,#1) \, dx}
\newcommand{\N}{{\mathbb{N}}}

\newcommand{\R}{{\mathbb{R}}}

\newcommand{\Ecal}{{\mathcal{E}}}

\newcommand{\Kcal}{{\mathcal{K}}}

\DeclareMathOperator{\dist}{dist}
\DeclareMathOperator{\supp}{supp}

\DeclareMathOperator{\dive}{div}

\title{Ergodic Mean Field Games: existence of local minimizers up to the 
Sobolev critical case}
\author{Marco Cirant, Alessandro Cosenza, Gianmaria Verzini}

\begin{document}
\maketitle

\begin{abstract} 
We investigate the existence of solutions to viscous ergodic 
Mean Field Games systems in bounded domains with Neumann boundary conditions 
and local, possibly aggregative couplings. In particular we exploit the 
associated variational structure and search for constrained minimizers of a 
suitable functional. Depending on the growth of the coupling, we detect the existence of global minimizers in the mass subcritical and critical case, 
and of local minimizers in the mass supercritical case, notably up to the Sobolev critical case.
\end{abstract}
\noindent
{\footnotesize \textbf{AMS-Subject Classification}}. 
{\footnotesize 35J50, 35Q89, 35B33, 49N80}

\section{Introduction} 

In this work we investigate the existence of solutions to the following system arising in the theory of viscous ergodic Mean Field Games, with Neumann boundary conditions and local coupling
\begin{equation}
\label{e:sys}
\begin{cases}
-\Delta u +H(\nabla u )+ \lambda = f(x,m(x)) \qquad & \text{on} \ \Omega \\
-\Delta m - \dive({m\nabla H (\nabla u )})=0 \qquad  & \text{on} \ \Omega \smallskip\\
\dfrac{\partial u }{\partial n}=0,  \ \  
\dfrac{\partial m}{\partial n } + m \nabla H(\nabla u )\cdot n =0 \qquad & \text{on} \ \partial \Omega
\smallskip\\
\int_{\Omega}m=1, \qquad \int_{\Omega} u=0,
\end{cases}
\end{equation}
and their minimality properties in a suitable variational framework. Throughout the paper, $\Omega$ is a bounded domain of $\mathbb{R}^{N}$, $N \ge 1$, with boundary of class ${C}^{3}$. To simplify some computations, we suppose that $|\Omega| = 1$, though all the arguments work for  $|\Omega| \neq 1$. On the Hamiltonian $H$, we assume that
\begin{flalign}&
    \begin{aligned}
    \label{e:ass:H}
    & C_{H}^{-1}\abs{p}^\gamma-K_{H}\leq H(p)\leq C_{H}\abs{p}^\gamma+K_H\\
    & \nabla H(p)\cdot p -H(p)\geq C^{-1}_{H}\abs{p}^\gamma-C_{H} \\
    & \abs{\nabla H(p)}\leq C_{H}\abs{p}^{\gamma -1}+K_H
    \end{aligned}
    \end{flalign}
for some $C_{H}>1,K_{H}>0,\gamma>1$. 
As for the coupling, we suppose that $f:\mathbb{R}^{N}\times[0,+\infty) \longrightarrow\mathbb{R}$ is locally Lipschitz and 
\begin{equation}
\label{e:ass:f}
    -\cf m^{q-1}-K_{f} \leq f(x,m)\leq \cf m^{q-1}+K_{f}
\end{equation}
for some $\cf>0$, $K_{f}>0$, $q>1$. 

\medskip

Mean-Field Games (MFG) have been introduced in the seminal papers by Lasry and Lions \cite{lasry_lions_2007_II} and Huang, Caines and Malham\`e, \cite{huang_malhame_caines_2006}, with the aim of describing Nash equilibria in differential games with infinitely many indistinguishable agents. The system \eqref{e:sys} characterizes these equilibria in an ergodic game, where the cost of a typical agent is averaged over an infinite-time horizon. Neumann boundary conditions come from the assumption that agents' trajectories are constrained to $\Omega$ by normal reflection at the boundary (as in \cite{cir_goffi_leonori,cirant_verzini_2017}).

\medskip

Though systems of type \eqref{e:sys} have been widely investigated over the last decade \cite{bardi_feleqi, bernardini_cesaroni, cesa_cirant_2019, cirant_2014, cirant_2016,  GPV, pime_voska}, it is not yet known the existence of classical solutions in the full range $q > 1$. While known restrictions on $q$ may be artificial, they are sometimes \textit{structural}. Let us briefly describe the purely quadratic case $H(p) = |p|^2$ to clarify this point. By the classical Hopf-Cole transformation $\phi=e^{-u} / \int e^{-u} = \sqrt m$, \eqref{e:sys} boils down to
\[
\label{e:sysq}
\begin{cases}
-\Delta \phi = \lambda \phi - f(x,\phi^2)\phi \qquad & \text{on} \ \Omega \\
\frac{\partial \phi }{\partial n}=0 \qquad & \text{on} \  \partial \Omega \\
\int_{\Omega}\phi^2\mathrm{d} x=1,
\end{cases}
\]
which, for $F'=f$,  are the Euler-Lagrange equations of the functional defined on $W^{1,2}(\Omega)$
\[
\mathcal F(\phi)= \int_\Omega |\nabla \phi|^2 + F(\phi^2) \mathrm{d} x \qquad \text{constrained to $\int_{\Omega}\phi^2\mathrm{d} x=1$.}
\]
In the model case $F(m) = C_f m^{q}$, it is well known (by the Gagliardo-Nirenberg inequality) that $2q \le 2 + 4/N$ is necessary for $\mathcal F$ to be bounded from below when constrained to $\int \phi^2 = 1$, while $2q < 2 + 4/(N-2)$ is needed, in dimension $N\ge3$, for the compact embedding of $W^{1,2}(\Omega)$ into $L^{2q}(\Omega)$, which is crucial to set up variational methods. There are therefore two critical values that determine three regimes, each one of exhibiting different properties and difficulties (see also \cite{MR4191345} and references therein). Note that when $C_f > 0$ the previous functional is convex, but we are not assuming here this property. In fact, we will be mainly interested in a model case $C_f < 0$, that corresponds to a MFG where aggregation between agents is enforced.

In the general nonquadratic case, \eqref{e:sys} borrows the following variational structure (see for example \cite{cesaroni_cirant_2019, Santambrogio2020} and references therein for further details). Let
\begin{equation}
\label{e:enrg}
\mathcal{E}(m,w):=
\begin{cases}
\displaystyle \int_{\Omega} \ml{m}{w} + F(x,m) \,dx & \text{if} \ (m,w)\in \mathcal{K}\\
+\infty & \text{otherwise,} 
\end{cases}
\end{equation}
where
\begin{equation}
\label{e:F}
F(x,m):=
\begin{cases}
\displaystyle \int_{0}^{m}f(x,n) \ dn & \text{if} \ m\geq 0 \\
+\infty & \text{if} \ m< 0,
\end{cases}
\end{equation}
and
\begin{equation}
\label{e:const}
    \begin{split}
        \mathcal{K}:=\big\{&(w,m)\in L^{\frac{\gamma' q }{\gamma' +q -1}}(\Omega)\cap {W}^{1,r}(\Omega)\times\leb{1} \ \text{s.t.}\big. \\ &\int_{\Omega}\nabla m \cdot \nabla \phi \, dx =\int_{\Omega}w\cdot \nabla\phi \, dx \ \  \forall \phi \in C^{\infty}(\overline{\Omega}),\\&\big. \int_{\Omega}m \,dx \ =1, m\geq 0 \ \text{a.e.} \big\}, \qquad \text{where } \  \frac{1}{r}:=\frac{1}{\gamma'}+\frac{1}{\gamma q}.
    \end{split}
\end{equation}
It has been first observed in \cite{cirant_2016} that for any $\gamma > 1$, two critical values of $q$ can be identified:
\[
\Bar{q}=1 + \frac{\gamma'}{N} \text{ \ (mass critical)}, \qquad \qc = 
\begin{cases}
1+\frac{\gamma'}{N-\gamma'} & \text{if $\gamma' < N$} \\
+ \infty & \text{if $\gamma' \ge N$}
\end{cases}
 \text{ \ (Sobolev critical)}.
\]
As in the quadratic case, $q \le \bar q$ is necessary for $\mathcal{E}$ to be bounded from below, while $q < q_c = r^*$ guarantees the compact embedding of $W^{1,r}$ into $L^q$ (this explains why we call this exponent ``Sobolev critical''). 

Note that if $\gamma = 2$, $\bar q$ and $q_c$ agree with the critical exponents mentioned in the previous paragraph.

Most of the analysis on systems of type \eqref{e:sys} has been carried out in the case $q < q_c$. By means of fixed point methods, solutions have been shown to exist in \cite{cirant_2016} in the range $\bar q \le q < q_c$ under the further assumption that $C_f$ be small enough (for problems which are set on the flat torus). No results are known when $q = q_c$, and if $q > q_c$ solutions may even fail to exist, at least when the system is set on the whole euclidean space.

The main goal of the present paper is to show the existence of solutions that are \textit{local} minimizers of $\mathcal{E}$, in particular in the range $\bar q \le q \le q_c$, up to the critical exponent $q_c$ The further property of (local) minimality obtained here may be an important feature in the study of their \textit{stability}, that is, their ability to capture the long-time behavior of the parabolic version of \eqref{e:sys}.
Our main result reads as follows.
\begin{theorem}
\label{th:existence}
Assume that \eqref{e:ass:H} and \eqref{e:ass:f} hold.
Suppose that either 
\begin{enumerate}
    \item $1<q<\Bar{q}$, or
    \item $q=\Bar{q}$ and  $\cf \le C_{\bar{q}}$, or
    \item $\Bar{q}< q< \qc$ and $\cf \le C_{q}$, or
    \item $q=\qc$ and $\cf, K_f, K_H \le C_{crit}$.
\end{enumerate}
Then there exists a solution $(u,\lambda,m)$ to the system \eqref{e:sys}, and $m>0$.

Moreover, the pair $(m,-m\nabla H(\nabla u))$ is a global minimizer of $\Ecal$ on $\Kcal$ in cases 1 and 2, or a minimizer of $\Ecal$ on $\mathcal{K}\cap \{m:\|u\|_{L^q}^q<\bar\alpha\}$, for a suitable $\bar\alpha>0$, in cases 3 and 4.
\end{theorem}

A solution to the system \eqref{e:sys} is a triple $(u,\lambda,m)\in {C}^{2,\theta}(\overline{\Omega})\times \mathbb{R}\times {W}^{1,p}(\Omega)$ for all $\theta \in (0,1)$ and all $p>1$ such that $(u,\lambda)$ is a classical solution to the Hamilton-Jacobi equation and $m$ is a weak solution to the Fokker-Planck equation. The constants $C_{\bar q}, C_{q}$ and $C_{crit}$ are explicitly calculated throughout Section \ref{sec:existence}. We mention here that they depend on the data $q,C_{H},\gamma$ and also on some regularity and embedding constants $C_{S},C_{E},C_{q},\delta_{0}$, which in turn depend on $q, \Omega, N$. 

The model Hamiltonian $H(p) = |p|^\gamma$, $\gamma >1$ clearly falls into our set of assumptions, and we can allow for a general growth of type $|p|^\gamma$ with different coefficients from above and below. We just need to be careful when $q = \qc$: $H$ has to be small  enough for small values of $|p|$. As for $f$, we include the model case
\[
f(x,m)=\cf a(x) m^{q-1} + K_f b(x),
\]
where $a,b$ are smooth functions (with no sign condition). When $\bar q \le q < \qc$ we just need $C_f$ to be small enough, while $q=\qc$ requires further smallness conditions on $K_f$.

\smallskip

The way local minima are identified is inspired by \cite{MR3689156,noris_tavares_verzini}. The functional $\Ecal$ is minimized first on the intersection between the constraint $\mathcal K$ and the ball $\{m:\|u\|_{L^q}^q \le \bar\alpha\}$. Estimates based on the Gagliardo-Nirenberg inequality for competitors belonging to $\mathcal K$, which involves a differential constraint of Fokker-Planck type, allow to choose $\bar \alpha$ in a way that these minimizers lie in fact in the open ball $\{m:\|u\|_{L^q}^q < \bar\alpha\}$, and therefore they give rise to solutions to the optimality conditions \eqref{e:sys}. There are a few technical obstacles to produce a solution of \eqref{e:sys} from a local minimizer of $\Ecal$. These are worked out following a strategy that involve a regularization of the functional and a linearization which allows to use convex duality methods. This strategy is detailed for example in \cite{cesaroni_cirant_2019}.

Note that he critical case requires further smallness assumptions on the coefficients. There is an interesting connection between this endpoint case and another phenomenon of criticality arising in Hamilton-Jacobi equations. If $m$ is in a ball of $L^q$, then $f(m)$, which is the right-hand side of the Hamilton-Jacobi equation in \eqref{e:sys}, will be bounded in $L^{\gamma'/N}$. Recent works \cite{cirant_goffi_2021, cir_verzini_2022} on the so-called maximal regularity of Hamilton-Jacobi equations underlined the crucial role of the exponent $\gamma'/N$: if the right-hand side $f(m)$ is uncontrolled in $L^{\gamma'/N}$, then there is no way of controlling $H(\nabla u), \Delta u$ separately in $L^{\gamma'/N}$, and therefore there is no hope to deduce further regularity via bootstrap arguments. In this sense, $\gamma'/N$ is critical. If the $L^{\gamma'/N}$ norm of $f(m)$ is \textit{small enough}, then $H(\nabla u), \Delta u$ can be controlled in the same Lebesgue space, which leads to further regularity. Additional smallness assumptions should be then expected when dealing with the endpoint case $q=\qc$.

\medskip

A natural question concerns the uniqueness of solutions to \eqref{e:sys}, which is not expected in general. We discuss here a few examples, limiting to the model nonlinearities
\[
H(p)=\frac{1}{\gamma}\abs{p}^{\gamma},\qquad
L(q)=\frac{1}{\gamma'}\abs{q}^{\gamma'},\qquad
f(x,m) = f(m) = \pm \cf m^{q-1}.
\]
In this setting, clearly the system \eqref{e:sys} admits the trivial solution $(u_{tr},\lambda_{tr},m_{tr})\equiv (0,f(1),1)$, which has a corresponding energy $\mathcal{E}(m_{tr},-m_{tr}\nabla H(\nabla u_{tr}))=\mathcal{E}(1,0)=F(1)$.

First of all, if $f(m) = \cf m^{q-1}$, $\cf>0$, it is well known that \eqref{e:sys} admits 
a unique solution, which must coincide with the trivial one (see for example \cite{lasry_lions_2007_II}). Instead, let us  
assume $f(m) = -\cf m^{q-1}$, $\cf>0$: then existence of multiple solutions follows whenever we can build 
an appropriate couple $(m,w)$ in $\Kcal$ (or in $\Kcal\cap B_{\bar{\alpha}}$) satisfying 
$\mathcal{E}(m,w)<F(1)$; indeed, in such case the solution found in Theorem \ref{th:existence} cannot 
be the trivial one. To this aim, take any  $\phi \in C^2(\overline{\Omega})$ such that $\int_{\Omega}
\phi\, dx=0$, and consider $m=1+\eps \phi$ and $w=\nabla m$. For $0\le \eps\le\bar\eps$ small enough, 
we have that $(m,w)\in K\cap B_{\bar{\alpha}}$ and $f '(1-\bar\eps\|\phi\|_\infty) \le -C\cf<0$.
We obtain, for some $0<\xi<\bar\eps$,
\[ 
\begin{split}
    \mathcal{E}(m,w)&=\int_{\Omega}\ml{m}{w}\,dx+\int_{\Omega}F(m)\,dx\\
   	&= \frac{1}{\gamma'}\int_{\Omega}(1+\eps\phi(x))^{1-\gamma'}\abs{\eps \nabla \phi (x)}^{\gamma'}\,dx+\int_{\Omega}F(1+\eps \phi(x))\,dx \\
	&\leq 
    C\eps^{\gamma'}+F(1)+\eps f(1)\int_{\Omega}\phi(x)\,dx+\frac{\eps^{2}}{2}\int_{\Omega}
    f'(1+\xi\phi(x))\phi^{2}(x)\,dx \\ &\leq  F(1)+C\eps^{\gamma'}-C'\cf\eps^{2}.
    \end{split}
\] 
Now, this last quantity is strictly smaller than $F(1)$ either for $\gamma'>2$ and $\eps$ small 
enough, or for $\gamma'\le2$, $\cf > C\bar\eps^{-(2-\gamma')}/C'$ and $\eps = \bar\eps$. We deduce 
multiplicity of solutions when
\[
\text{either }\gamma<2,\qquad \text{ or } \gamma\ge2,\ 1<q<\bar q\text{ and $\cf$ sufficiently large} 
\]
(in case $\gamma\ge2$ and $q\ge\bar q$ one has to compare the conditions on $\cf$ here and in Theorem 
\ref{th:existence}). On the other hand, if $\gamma\ge2$ and $\cf$ is small enough, then uniqueness may be 
expected, see \cite{cirant_porretta}.

In case of multiplicity of solutions, a further question concerns the uniqueness of the minimizers. 
This is an interesting question, which will be the object of subsequent studies.

\bigskip

\textbf{Notations.}
For $k\in \N$ and $p\geq 1$ we denote by $\norm{u}_{p}$ and $\norm{u}_{k,p}$ the usual $\leb{p}$ and $W^{k,p}(\Omega)$ norm respectively. For $p\geq 1$, the exponent $p'$ is the conjugate exponent of $p$, $p'=\frac{p}{p-1}$. $C,C'$ and so on denote non-negative universal constants, which we need 
not to specify, and which may vary from line to line.

\section{Preliminaries}

\subsection{The Lagrangian}
\label{subsec:prelim}
The Legendre-Fenchel transform $L_H$ of $H$ is necessary for the construction of the energy associated to the system: 
\[
    L=L_{H}(q):=\sup_{p \in \mathbb{R}^{N}}[p\cdot q- H(p)]
\]
Under our assumptions on $H$, the following properties of $L_{H}$ are standard:
\begin{proposition}
\label{prop:leg:prop}
There exists $C_{L}>0$ such that for all $p,b \in \mathbb{R}^{N}$
\begin{enumerate}
    \item $L_{H}\in C^2\left(\mathbb{R}^{N}\backslash \{ 0 \} \right) $ and it is strictly convex.
    \item $0\leq C_{L}\abs{q}^{\gamma'}\leq L_{H}(q)\leq C_{L}^{-1}(\abs{q}^{\gamma'}+1)$
    \item $\nabla L_{H}(q)\cdot q - L_{H}(q)\geq C_{L}\abs{q}^{\gamma'}-C_{L}^{-1} $
    \item $C_{L}\abs{q}^{\gamma ' -1}-C_{L}^{-1}\leq \abs{\nabla L_{H}(q)}\leq  C_{L}^{-1}(\abs{q}^{\gamma ' -1}+1)$
\end{enumerate}
\end{proposition}
\begin{proof}
See for example \cite{cirant_2014_II}.
\end{proof}

The energy functional $\mathcal E$ involves the following Lagrangian term.

\begin{proposition}
\label{prop:mLF:prop}
The function
\begin{equation*}
    (m,w)\longrightarrow \ml{m}{w} = \begin{cases}
mL(-\frac{w}{m}) & \text{if} \ m>0 \\
0 & \text{if}\ m=0,w=0,\\
+\infty & \text{otherwise,}
\end{cases}
\end{equation*}
is convex, and strictly convex if restricted to $m>0$. We have,
\begin{equation}
\label{e:mL:car}
mH(p)=
\begin{cases}
 \displaystyle\sup_{w\in \mathbb{R}^{N}} \left[-p\cdot w -\ml{m}{w}\right] & \text{if} \ m\neq 0 \\
0 & \text{if} \ m = 0.
\end{cases}
\end{equation}
Moreover
\begin{equation}
    \label{e:mL:prop}
    C_{L}\frac{\abs{w}^{\gamma'}}{m^{\gamma'-1}}\leq \ml{m}{w}\leq C_{L}^{-1}\frac{\abs{w}^{\gamma'}}{m^{\gamma'-1}}+C_{L}^{-1}m.
\end{equation}
\end{proposition}
\begin{proof}
Equation \eqref{e:mL:car} is standard, see for example \cite{cesaroni_cirant_2019}.
Estimate \eqref{e:mL:prop} comes directly from Proposition \ref{prop:leg:prop}.
\end{proof}

\subsection{Fokker-Planck equations}
\label{sec:statfp}

We deal here with the Fokker-Planck equation 
\begin{equation}
\label{e:fp}
    \begin{cases}
    -\Delta m - \dive({m b})=0 \qquad  & \text{on} \ \Omega \\
     \frac{\partial m}{\partial n } + m  b\cdot n =0 \qquad & \text{on} \ \partial \Omega\\
    \int_{\Omega}m=1
    \end{cases}
\end{equation}
where $b:\mathbb{R}\rightarrow\rn$ will be (at least) in $L^{s}(\Omega;\rn)$, for some $s > N$. Solutions $m\in W^{1,2}(\Omega)$ will be in the standard weak sense: 
\begin{equation}
\label{e:fpsolution}
     \int_{\Omega}\nabla m \cdot \nabla \phi \, dx =\int_{\Omega}bm\cdot \nabla \phi \, dx \qquad \forall \phi \in W^{1,2}(\Omega),
\end{equation}
with $\int_{\Omega}m=1$. The following existence result for $b$ in $L^\infty$ is classical.

\begin{theorem}
\label{th:fp}
Let $b\in L^{\infty}(\Omega;\rn)$. Then there exist an unique weak solution $m$ to \eqref{e:fp}, $m \in W^{1,p}(\Omega)$ for every $p$ and
\begin{equation*}
    \norm{m}_{1,p}\leq C= C(\norm{b}_{\infty}, p, N, \Omega).
\end{equation*}
Moreover, $m\in C^{0,\alpha}(\overline{\Omega})$
for all $\alpha \in (0,1)$ and there exists $c=c(\norm{b}_{\infty}, p, N, \Omega)>0$ such that
\begin{equation*}
    c^{-1}\leq m(x) \leq c 
\end{equation*}
for all $x\in \Omega$.
\end{theorem}

\begin{proof}
see \cite[Th.\ II.4.4, II.4.5, II.4.7]{bensoussan_1988}.
\end{proof}

We investigate further regularity properties, and recall the following proposition, which is an useful $W^{1,p}$ regularity result for linear equations in divergence form.
\begin{proposition}
\label{prop:lpest}
Let $\rho \in \leb{p}$, with $p>1$.
Suppose that
\begin{equation}
    \label{e:lpest:hyp}
    \abs{\int_{\Omega}\rho\Delta\phi\, dx }\leq K\norm{\nabla \phi}_{p'}
\end{equation}
for all $\phi \in C^{\infty}(\overline{\Omega}), \frac{\partial \phi }{\partial n}=0$ for some $K>0$. Then $\rho\in W^{1,p}(\Omega)$ and there exists $C_{E} = C_E(N,\Omega,p)$ such that
\begin{equation}
    \label{e:lpest}
    \norm{\rho}_{1,p}\leq C_{E}(K+\norm{\rho}_{p}).
\end{equation}
Moreover, the same estimate holds in a local form, that is, for every $B_{2R} \subset \Omega$,
\begin{equation}
    \label{e:lpestloc}
    \norm{\rho}_{W^{1,p}(B_R)}\leq C_{E,R}(K+\|\rho\|_{L^p(B_{2R})}).
\end{equation}
\end{proposition}
\begin{proof}
See \cite[Theorems 7.1 and 8.1]{agmon_1959}. 
\end{proof}

As a straightforward consequence, given $w\in L^{p}(\Omega;\rn)$, $p>1$, any weak solution $m\in W^{1,p}(\Omega)$ to 
\begin{equation}
\label{e:fpprelim}
     \begin{cases}
    -\Delta m - \dive w=0 \qquad  & \text{on} \ \Omega \\
    \frac{\partial m}{\partial n } + w\cdot n =0 \qquad & \text{on} \ \partial\Omega\\
   \int_{\Omega}m=1.
    \end{cases}
\end{equation}
satisfies
\begin{equation}\label{e:cebound}
    \norm{m}_{1,p}\leq C_{E}(\norm{w}_{p}+\norm{m}_{p}).
\end{equation}

The previous estimate, combined with the Gagliardo-Nirenberg inequality, yields the following crucial result. 

\begin{proposition}
\label{prop:est}
Let $(m,w)\in L^1 \cap W^{1,r}(\Omega) \times L^1(\Omega)$ be a solution of \eqref{e:fpprelim}, and $E$ be defined by
\begin{equation*}
    \label{e:E}
    E:=\int_{\Omega}\frac{\abs{w}^{\gamma'}}{m^{\gamma'-1}}\, dx,
\end{equation*}
which is assumed to be finite. Then there exists $\cpr>0$ such that
\begin{align}
\label{e:estw}
    &\norm{m}_{1,r}\leq \cpr(E+1)\\
\label{e:estq}
    &\norm{m}_{q}\leq \cpr(E+1).
\end{align}
Moreover, if $q<\bar{q}$, there exists also $\delta >0$ such that
\begin{equation}
    \label{e:estq2}
    \norm{m}_{q}^{q(1+\delta)}\leq \cpr(E+1).
\end{equation}
Finally, if $q=\Bar{q}$, then \eqref{e:estq2} holds with $\delta=0$.

\end{proposition}
\begin{proof} The proof can be found for example in \cite{cesaroni_cirant_2019}, and its adaption to the problem with Neumann conditions is straightforward.
\end{proof}

To conclude the section, we present a sharper estimate which will be useful in an endpoint case of our analysis, and requires $b$ to be controlled in $L^{N}(\Omega)$ only. 
\begin{proposition}
\label{prop:fpestcrit}
Let $b\in L^{\infty}(\Omega;\rn)$ such that $\norm{b}_{N}\le \frac{1}{2C_{E}C_{S}} $, where $C_{E}$ is defined in \eqref{e:lpest} and $C_{S}$ is the $p$ Sobolev Embedding constant. Then for any $p<N$ the solution $m$ of \eqref{e:fp} satisfies.
\begin{equation}
    \label{e:fpestcrit}
    \norm{m}_{1,p}\leq C = C(p, N, \Omega)
\end{equation}
\end{proposition}
\begin{proof}
By \eqref{e:cebound} we have that
\begin{equation*}
    \norm{m}_{1,p}\leq C_{E}(\norm{bm}_{p}+\norm{m}_{p}).
\end{equation*}
By  H\"oder inequality, Sobolev Embedding and Interpolation
\begin{equation*}
     \norm{bm}_{p}+\norm{m}_{p}\leq \norm{b}_{N}\norm{m}_{p^{*}}+\norm{m}^{1-\theta}_{1}\norm{m}_{p^{*}}^{\theta}\leq C_{S}\norm{b}_{N}\norm{m}_{1,p}+\norm{m}_{1,p}^{\theta},
\end{equation*}
for some $\theta=\theta(p,N) < 1$, from which we can deduce the thesis using the bound assumed on $\norm{b}_{N}$.
\end{proof}

\subsection{Hamilton-Jacobi equations}
\label{subsec:stathj}
We now consider the Hamilton-Jacobi equation
\begin{equation}
\label{e:hj}
  \begin{cases}
  -\Delta u +H(\nabla u )+ \lambda = f(x) \qquad & \text{on} \ \Omega \\
\frac{\partial u }{\partial n}=0 \qquad & \text{on} \  \partial \Omega \\
\int_{\Omega}u=0 .
  \end{cases}  
\end{equation}
A solution of \eqref{e:hj} is a pair $(u,\lambda) \in C^{2}(\overline{\Omega})\times \R$ that satisfies the equation pointwise. 
\begin{theorem}
\label{th:hjexistence}
Let $f\in C^{\alpha}(\overline{\Omega})$, $0< \alpha < 1$. Then there exists an unique constant $\lambda \in \mathbb{R}$ such that \eqref{e:hj} has a unique solution in $C^{2,\alpha}(\overline{\Omega})$ and
\begin{equation}
    \label{e:deflambda}
    \lambda=\sup \left\{c\in \mathbb{R} \text{ s.t. } \exists u\in C^{2}(\Bar{\Omega}), \frac{\partial u}{\partial n}=0 \text{ on } \partial\Omega : -\Delta u +H(\nabla u) +c \leq f\right\}.
\end{equation}
Moreover, the following estimates hold:
\begin{align*}
 \norm{\nabla u}_{\infty}\leq&  K_1 = K_1(N, \Omega, \|f\|_\infty), \\
 \|u\|_{C^{2,\alpha}(\overline \Omega)} \leq & K_2 = K_2(N, \Omega, \alpha, \|f\|_{C^\alpha(\overline \Omega)}).
\end{align*}
\end{theorem}
Note that the estimate holds also locally, that is, for every $B_{2R} \subset \Omega$,
\[
\norm{\nabla u}_{L^\infty(B_R)}\leq  K_1 = K_1(N, \Omega, R, \|f\|_{L^\infty(B_{2R})}).
\]
The proof of this theorem is well-known in ergodic control theory, and it is typically obtained via a limiting procedure involving a discounted problem. The crucial gradient estimate which allows to pass to the limit in the procedure can be derived using the Bernstein method, see for example  \cite{cesaroni_cirant_2019, cirant_2014, lions_1985} and references therein. Though we are not going to use directly the characterization \eqref{e:deflambda} of $\lambda$ here, it is in fact a key step in the existence argument of Theorem \ref{th:duality}, which follows some standard lines involving convex duality.

We are now interested in finding additional regularity results, which will be used in our critical endpoint case. In particular, we are interested in finding bounds for $|\nabla u|^\gamma$ in $L^{{N}/{\gamma'}}$ depending on $f$ in $L^{{N}/{\gamma'}}$. This is a delicate endpoint case of the so-called $L^q$-maximal regularity for Hamilton-Jacobi equations, which has been recently discussed in \cite{cir_goffi_leonori, cir_verzini_2022} for $q \ge \frac{N}{\gamma'}$. We provide here a simple proof which exploits a smallness condition on $f$.

We start by introducing the following lemma for linear equations, and provide its standard proof for the reader's convenience.
\begin{lemma}
\label{lem:hjlemma}
Let $f\in L^{p}(\Omega)$ for $p>1$ and let $(u,\lambda)\in W^{2,p}(\Omega)\times \R$ be a solution of
\begin{equation}
\label{e:hjlemma}
    \begin{cases}
    -\Delta u +\lambda=f \text{ on } \Omega \\
    \frac{\partial u}{\partial n}=0  \text{ on } \partial \Omega\\
    \int_{\Omega} u=0.
    \end{cases}
\end{equation}
Then there exists $C=C(N,\Omega,p)>0$ (independent of $\lambda$) such that
\begin{equation}
    \norm{u}_{2,p}\leq C \norm{f}_{p}.
\end{equation}
\end{lemma}
\begin{proof}
Using elliptic regularity (which holds up to the boundary by homogeneous Neumann boundary conditions), we have that 
\begin{equation}
\label{e:hjlemmaregularity}
    \norm{u}_{2,p}\leq C(\norm{f-\lambda}_{p}+\norm{u}_{p}).
\end{equation} 
If we test \eqref{e:hjlemma} with a constant function, we get that necessarily 
\begin{equation}
    \lambda =\frac{1}{\abs{\Omega}}\int_{\Omega}f\,\mathrm{d} x,
\end{equation}
hence we can suppose that $\norm{f-\lambda}_{p}\leq C\norm{f}_{p}$.
Let us now claim that 
\begin{equation}
    \norm{u}_{p}\leq C \norm{f}_{p}.
\end{equation}
Indeed, suppose by contradiction that there exists a sequence $(u_{n},f_{n})$ satisfying \eqref{e:hjlemma} such that 
\begin{equation}
    \norm{u_{n}}_{p}>n\norm{f_{n}}_{p}.
\end{equation}
By \eqref{e:hjlemmaregularity} we get 
\begin{equation}
    n\norm{f_{n}}_{p}\leq C(\norm{f_{n}}_{p}+\norm{u_n}_{p}),
\end{equation}
hence we can conclude that $\frac{\norm{f_{n}}_{p}}{\norm{u_{n}}_{p}}\rightarrow 0$. Now let $v_{n}=\frac{u_{n}}{\norm{u_{n}}_{p}}$. Clearly $\norm{v_{n}}_{p}=1$, $\int_{\Omega}v_{n}=0$ and 
\begin{equation*}
    \norm{v_{n}}_{2,p}\leq C\left(\frac{\norm{f_{n}}_{p}}{\norm{u_{n}}_{p}}+1\right)\leq C.
\end{equation*}
Hence we have that there exists a subsequence $v_{n}\rightarrow v$ strongly in $L^{p}(\Omega)$ and weakly in $W^{2,p}(\Omega)$. Moreover, $\norm{v}_{p}=1$, $\int_{\Omega}v=0$. Now since $u_{n}$ satisfies \eqref{e:hjlemma} we have for all $\phi\in C^{\infty}(\overline{\Omega})$
\begin{equation*}
    \int_{\Omega}\nabla v_{n}\nabla \phi\,dx=\frac{1}{\norm{u_{n}}_{p}}\int_{\Omega}(f_{n}-\lambda)\phi\,dx\leq \frac{C \norm{f_{n}}_{p}}{\norm{u_{n}}_{p}} \norm{\phi}_{p'}.
\end{equation*}
Passing to the limit we get that $ \int_{\Omega}\nabla v\nabla \phi\,dx=0$ for all $\phi\in C^{\infty}(\overline{\Omega})$, hence $\nabla v=0$ which implies $v=K$ with $K=0$ by the constraint. This is a contradiction with $\norm{v}_{p}=1$.
\end{proof}
Using this estimate the idea is to construct a "barrier" for $\nabla u $ which we employ in a topological fixed point argument. This means finding a value $M$ such that no solutions exists with $\norm{\nabla u}_{N(\gamma'-1)}=M$.
\begin{proposition}
\label{prop:hjbarrier}
Under the hypothesis of Theorem \ref{th:hjexistence}, let $(u,\lambda)$ be the unique solution of \eqref{e:hj}. Then there exists $\delta_{0}$ depending on $N, \Omega, \gamma, p, H$ such that if \[K_H+\norm{f}_{\frac{N}{\gamma'}}\le\delta \le \delta_{0},\] then
\begin{equation}
\label{e:hjbarrier}
    \norm{\nabla u }_{N(\gamma-1)}< M \quad \text{or} \quad  \norm{\nabla u }_{N(\gamma-1)}> M
\end{equation}
for some $M$ depending on $\delta$ (and $N, \Omega, \gamma, p, C_H$). Moreover, $M=M(\delta) \to 0$ as $\delta \to 0$.
\end{proposition}
\begin{proof}
Let $p=\frac{N}{\gamma'}$, so $p^{*}=p \gamma$. Hence by Sobolev embedding, Lemma \ref{lem:hjlemma} and assumptions on $H$,
\begin{equation}
     y:=\norm{\nabla u}_{p\gamma}\leq \norm{u}_{1,p\gamma}\leq C \norm{u}_{2,p}
    \leq C(\norm{\nabla u}_{p\gamma}^{\gamma}+K_H +\norm{f}_{p})\leq C y^{\gamma}+C'\delta_{0}.
\end{equation}

We now choose $\delta_{0}$ such that there exists two positive solutions $y_{1},y_{2}$ to $y=Cy^{\gamma}+C'\delta_{0}$. In this way we can choose $M$ such that $y_{1}<M<y_{2}$. Moreover, $y_{1}$ can be made arbitrarily small if $\delta \le \delta_{0}$ is chosen small enough, hence also $M>0$ can be made arbitrarily small.
\end{proof}
Now that we have the barrier, we can use  Leray-Schauder fixed point theorem to deduce the following regularity result (which is in fact also an existence result). 
\begin{theorem}
\label{th:hjapriori2}
Under the hypothesis of Theorem \ref{th:hjexistence}, there exists $\delta_{0}>0$ such that if $K_H+\norm{f}_{\frac{N}{\gamma'}}\le \delta \le\delta_{0}$ then
\begin{equation}
\label{e:hjapriori2}
    \norm{\nabla u }_{N(\gamma-1)}< M.
\end{equation}
Moreover, $M=M(\delta) \to 0$ as $\delta \to 0$ ($M,\delta,\delta_0$ are as in the previous proposition). 
\end{theorem}
\begin{proof}
We use Leray-Schauder fixed point theorem (see for instance \cite[Thm.\ 4.3.4]{MR2816471}). Let 
$U=\{w\in C^{2}(\overline{\Omega}): \norm{\nabla u }_{N(\gamma-1)}< M \}$ and consider the operator $T:\overline{U}\rightarrow C^{2}(\overline{\Omega})$ where $u=T(w)$ is defined by the solution of the system 
\begin{equation*}
  \begin{cases}
  -\Delta u +H(\nabla w )+ \lambda = f(x) \qquad & \text{on} \ \Omega \\
\frac{\partial u }{\partial n}=0 \qquad & \text{on} \  \partial \Omega \\
\int_{\Omega}u=0 .
  \end{cases}  
\end{equation*}

We claim that $T$ is continuous and compact. Under this condition, by homotopy with the identity, 
Leray-Schauder theorem asserts that either $T$ has a fixed point or there exists $s\in (0,1)$ and $u\in \partial \Omega$ such that $u=s T(u)$. However, the latter possibility is excluded by Proposition \ref{prop:hjbarrier} (replacing $H(\nabla u), f$ with $sH(\nabla u), sf$ respectively). Hence $T$ has a fixed point in $U$, which by Theorem \ref{th:hjexistence} is the unique solution of \eqref{e:hj}, from which we can deduce the desired estimate. 

Thus we are left to prove that $T$ is continuous and compact. 
Clearly $T_{1}:C^{2}(\overline{\Omega})\rightarrow C^{\alpha}(\overline{\Omega})$ defined by $u=H(\nabla u)$ is continuous by our assumptions on $H$. Moreover let $T_{2}:C^{\alpha}(\overline{\Omega})\rightarrow C^{2,\alpha}(\overline{\Omega})$ be the operator defined by $u=T_{2}z$ solving
\begin{equation*}
  \begin{cases}
  -\Delta u +z+ \lambda = f \qquad & \text{on} \ \Omega \\
  \frac{\partial u }{\partial n}=0 \qquad & \text{on} \  \partial \Omega \\
   \int_{\Omega}u=0 .
  \end{cases}  
  \end{equation*}
This is well defined and continuous by elliptic regularity and because we supposed $f\in C^{\alpha}$. Finally, by Ascoli-Arzel\`a the immersion $T_{3}:  C^{2,\alpha}(\overline{\Omega})\rightarrow C^{2}(\overline{\Omega})$ is continuous and compact. Hence we conclude that $T=T_{3}\circ T_{2}\circ T_{1}$ is continuous and compact.
\end{proof}

\section{Existence of a solution}
\label{sec:existence}

We are now ready prove Theorem \ref{th:existence}. 
The proof of this theorem consists of several steps, which are explained in detail in the following sections. Since we are dealing with local couplings, the duality procedure to get a solution to the system \eqref{e:sys} is rather delicate (as the dual problem to the minimization of $\mathcal E$ would be related to solutions of an Hamilton-Jacobi equation with rough right-hand side). We will introduce a family of regularized problem with smoothing couplings, associate their energies $\mathcal{E}_{\ep}$ and then prove the existence of (local) minimizers. Once we have a minimum point for the regularized energy, we will perform the convex duality argument to deduce the existence of a regular solution. The solution to the initial problem will be then found by an appropriate limit procedure on the regularized sequence of solutions.

\subsection{Regularization}
Let us consider, for $\ep>0$, the following regularized system
\begin{equation}
    \label{e:sysreg}
\begin{cases}
-\Delta u +H(\nabla u )+ \lambda = f_{\ep}[m](x) \qquad & \text{on} \ \Omega \\
-\Delta m - \dive({m\nabla H (\nabla u )})=0 \qquad  & \text{on} \ \Omega \\
\frac{\partial u }{\partial n}=0 \qquad & \text{on} \  \partial \Omega \\
\frac{\partial m}{\partial n } + m \nabla H(\nabla u )\cdot n =0 \qquad & \text{on} \ \partial \Omega\\
\int_{\Omega}m=1, \qquad \int_{\Omega}u=0, 
\end{cases}
\end{equation}
where
\begin{equation}
         \label{e:feps}
    f_{\ep}[m](x):=f\left(\cdot,m*\chi_{\ep}(\cdot)\right)*\chi_{\ep}(x) = \int_{\rn}\chi_{\ep}(x-y)f\left(y,\int_{\rn}
    m(z)\chi_{\ep}(y-z)\, dz \right)\, dy
\end{equation}
and $\chi_{\ep}$ is a sequence of standard symmetric mollifiers approximating the unit ($f$ and $m$ are extended to $0$ outside $\Omega$).
We notice that given 
\begin{equation}
    \label{e:Feps}
    F_{\ep}[m]:=\int_{\Omega}F(x,m*\chi_{\ep}(x))\, dx
\end{equation}
we have that it holds
\begin{equation}
    \label{e:l2der}
    F_{\ep}[m']-F_{\ep}[m]=\int_{0}^{1}\int_{\Omega}f_{\ep}[(1-t)m+tm'](x)(m'-m)(x)\, dx
\end{equation}
for $m,m'\in \leb{1}$ and $\int_{\Omega}m=\int_{\Omega}m'=1$. This means that the regularized problem also admits a potential. 
Let us observe that using the properties of mollifiers and the assumptions on $f$, 
 the following estimates hold:
\begin{equation}
\label{e:Feps:ass}
-\frac{\cf}{q}\norm{m}_{q}^{q}-K_{f}\leq F_{\ep}[m]\leq \frac{\cf}{q}\norm{m}_{q}^{q}+K_{f},
\end{equation}
and
\begin{equation}
\label{e:Feps:ceps}
-\frac{\cf}{q}\sup_{\Omega}\chi_{\ep}^{q}-K_{f}\sup_{\Omega}\chi_{\ep}\leq F(x,m*\chi_{\ep}(x))\leq \frac{\cf}{q}\sup_{\Omega}\chi_{\ep}^{q}+K_{f}\sup_{\Omega}\chi_{\ep}.
\end{equation}
We now introduce the energy of the approximated system:
\begin{equation}
\label{e:enrgeps}
\mathcal{E}_{\ep}(m,w):=
\begin{cases}
\displaystyle \int_{\Omega} \ml{m}{w}\, dx + F_{\ep}[m] & \text{if} \ (m,w)\in \mathcal{K}\\
+\infty & \text{otherwise}.
\end{cases}
\end{equation}

\subsection{Minimization of the regularized functional}
Our goal now is to find a minimizer for the energy of the regularized system. Here, the exponent $\Bar{q}$ comes into play. If $q< \Bar{q}$, the energy can be proven to admit a global minimum. This is the case addressed in \cite{cesaroni_cirant_2019}. If $q=\Bar{q}$, a global minimum can be found under some condition on the coefficients. If $q> \Bar{q}$, no global minima is present in general and we have to look for local minima. Following the idea introduced in \cite{MR3689156,noris_tavares_verzini}, let us introduce, for $\alpha\geq 1$
\begin{equation}
    B_{\alpha}:=\{  (m,w)\in \mathcal{K} : \norm{m}_{q}^{q}\leq \alpha\}
\end{equation}
and 
\begin{equation}
    U_{\alpha}:=\{  (m,w)\in \mathcal{K} : \norm{m}_{q}^{q}= \alpha\}.
\end{equation}
Let us also define, for $\eps>0$ fixed,
\begin{equation}
    c_{\alpha}=\inf_{(m,w)\in \mathcal{K}\cap B_{\alpha}}\mathcal{E}_{\ep}(m,w)
\end{equation}
and 
\begin{equation}
    \hat{c}_{\alpha}=\inf_{(m,w)\in \mathcal{K}\cap U_{\alpha}}\mathcal{E}_{\ep}(m,w).
\end{equation}
We start by proving the existence of a minimum in the sets $B_{\alpha}$, and also adapt the arguments to prove the existence of a global minimum under the stricter assumptions. 
\begin{lemma}
\label{le:minim}
For all $\alpha\geq 1 $, $c_{\alpha}$ is achieved. Moreover, if $q<\Bar{q}$ or $q=\Bar{q}$ and $\cf<C_{\Bar{q}}=q\frac{\cl}{\cpr}$ , then $\mathcal{E}_{\ep}$ has a global minimum on $\mathcal{K}$.
\end{lemma}
\begin{proof}
Let us begin by bounding $\mathcal{E}_{\ep}$ by below.
Using \eqref{e:mL:prop}, \eqref{e:estq} and \eqref{e:Feps:ass} we have:
\begin{equation}
\label{e:energybndal}
    \epsen{m}{w}\geq \frac{\cl}{\cpr}\norm{m}_{q}-\cl-\frac{\cf}{q}\norm{m}_{q}^{q}-K_{f}\geq K',
\end{equation}
since we know that $\norm{m}_{q}^{q}\leq \alpha$ in $B_{\alpha}$.
Consider now a minimizing sequence $(m_{n},w_{n})$ for $c_{\alpha}$. Eventually, $\epsen{m_{n}}{w_{n}}\leq c_{\alpha}+1$. Hence, again by 
\eqref{e:mL:prop} and \eqref{e:Feps:ass} we have
\begin{equation}
\label{e:E:estdimal}
    \int_{\Omega}\frac{\abs{w_{n}}^{\gamma'}}{m_{n}^{\gamma'-1}}\leq \cl^{-1}\left(c_{\alpha}+1-F_{\ep}[m_{n}]\right)\leq \cl^{-1}(c_{\alpha}+1+K_{f}+\frac{\cf}{q}\alpha)
\end{equation} 
which implies that $\left(\int_{\Omega}\frac{\abs{w_{n}}^{\gamma'}}{m_{n}^{\gamma'-1}}\right)$ is bounded.
Using \eqref{e:estw}, we have that $\norm{m_{n}}_{1,r}\leq C$. We can use Sobolev embeddings to conclude that up to subsequences
\begin{equation*}
m_{n}\rightarrow m \text{ a.e. on } \Omega, \quad
    m_{n}\rightarrow m \text{ on } \leb{1}, \quad 
    m_{n}\rightharpoonup m\text{ on }W^{1,r}(\Omega).
\end{equation*}
Then, using H\"older inequality
\begin{equation*}
    \int_{\Omega}\abs{w_{n}}^{\frac{\gamma' q }{\gamma' +q -1}}\, dx\leq \left(\int_{\Omega}\frac{\abs{w_{n}}^{\gamma'}}{m_{n}^{\gamma'-1}}\right)^{\frac{q}{\gamma'+q-1}}\norm{m_{n}}_{q}^{\frac{\gamma'-1}{q(\gamma'+q-1)}},
\end{equation*}
hence $w_{n}$ is equibounded in $\leb{\frac{\gamma' q }{\gamma' +q -1}}$ and so $w_{n}\rightharpoonup w$ in $\leb{\frac{\gamma' q }{\gamma' +q -1}}$. By $\leb{1}$ convergence of $m_{n}$ we can conclude in a standard way that $m \geq 0$ and that $\int_{\Omega}m=1$. Moreover, the convergences are strong enough to pass to the limit in the constraint $\mathcal{K}$, that is, $(m,w)\in \mathcal{K}$. Fatou's lemma also implies that $m \in B_{\alpha}$.
\\
To conclude, it is known that $\int_{\Omega}\ml{m}{w}\, dx$ is lower-semicontinuous with respect to the weak convergence of $W^{1,r}(\Omega)\times \leb{\frac{\gamma' q }{\gamma' +q -1}}$ (indeed, one can exploit its convexity and adapt classical results that connect convexity and lower semicontinuity, see for instance \cite[Th.\ 2.2.1]{evans_1990}). Moreover, using \eqref{e:Feps:ceps} and the Dominated Convergence Theorem we deduce that $F_{\ep}$ is strongly continuous with respect to the $\leb{1}$ convergence. 
Hence,
\begin{equation*}
    \epsen{m}{w}\leq \liminf_{n} \int_{\Omega}\ml{m_{n}}{w_{n}}\, dx + \lim_{n} F_{\ep}[m_{n}]  \leq \liminf_{n} \epsen{m_{n}}{w_{n}}=c_{\alpha}.
\end{equation*}
Now suppose that $q<\bar{q}$. Then \eqref{e:estq2} holds. The proof of the existence of a minimizer is completely analogous as before, but there is no need to restrict the set $B_{\alpha}$. Indeed, instead of \eqref{e:energybndal}, we can directly infer using \eqref{e:estq2} that 
\begin{equation}
    \label{e:energybnd}
     \epsen{m}{w}\geq \frac{\cl}{\cpr}\norm{m}_{q}^{q(1+\delta)}-\cl-\frac{\cf}{q}\norm{m}_{q}^{q}-K_{f}\geq K'.
\end{equation}
Moreover, we can set $e=\inf_{(m,w)\in \mathcal{K}}\epsen{m}{w}$ and argue as in \eqref{e:E:estdimal}, using \eqref{e:estq2} to conclude that
\begin{equation*}
    \label{e:E:estdim}
     \int_{\Omega}\frac{\abs{w_{n}}^{\gamma'}}{m_{n}^{\gamma'-1}}\leq \cl^{-1}\left(e+1+K_{f}+\frac{\cf}{q}\left(   \cpr \int_{\Omega}\frac{\abs{w_{n}}^{\gamma'}}{m_{n}^{\gamma'-1}}\,dx +1      \right)^{\frac{1}{1+\delta}}\right)
\end{equation*}
which again implies that $\int_{\Omega}\frac{\abs{w_{n}}^{\gamma'}}{m_{n}^{\gamma'-1}}$ is bounded.
Finally, if $q=\Bar{q}$ the previous steps are justified provided that $\cf<q\frac{\cl}{\cpr}$ and a global minimum exists. 
\end{proof}

\begin{remark}
Let $(m_{\ep},w_{\ep})$ be a minimizer $\mathcal{E}_{\ep}$ as in the previous lemma. Then, there exists $C>0$ independent of $\ep$
such that 
\begin{equation}
    \label{e:estqw}
    \norm{m_{\ep}}_{q}\leq C, \qquad \norm{m_{\ep}}_{1,r}\leq C 
\end{equation}
and 
\begin{equation}
\label{e:estqw2}
    \norm{w_{\ep}}_{\frac{\gamma'q}{\gamma'+q-1}}\leq C.
\end{equation}
Indeed, when $q\leq\Bar{q}$ we can use the fact that 
\begin{equation}
\label{e:estindip}
    \epsen{m_{\ep}}{w_{\ep}}\leq \epsen{1}{0}\leq \cl^{-1}+\frac{\cf}{q}+K_{f}
\end{equation}
to conclude that inequalities \eqref{e:estqw} and \eqref{e:estqw2} hold with a constant independent on $\eps$. If $q> \Bar{q}$, the same results hold, the proof being even easier since $m_{\ep}\in B_{\Bar{\alpha}}$, and $\Bar{\alpha}$ is independent of $\ep$.
\end{remark}

We see in the above lemma the role played by $\Bar{q}$. When $q>\Bar{q}$, $\mathcal{E}_{\ep}$ is not globally bounded from below and no global minima exist (see the remark below). To show that a local minimum exists in $B_{\alpha}$, we are left to prove that the candidate obtained in the previous lemma does not belong to $U_{\alpha}$. To this aim, we look for $\Bar{\alpha}>1$ such that $c_{\Bar{\alpha}}<\hat{c}_{\Bar{\alpha}}$.\\
\begin{remark}
\label{rem:umbond}

Let us notice that if $q\geq \bar{q}$, $\mathcal{E}$ (and also $\mathcal{E}_{\ep}$) could indeed be unbounded.

Let $f(x,m(x))=-\cf m(x)^{q-1}-K_{f}$ and $q>1+\frac{\gamma'}{N}$. Then, there exists $(m_{n},w_{n})$ such that $\mathcal{E}(m_{n},w_{n})\rightarrow -\infty$.
Choose $m_{0} \in C^{\infty}_{0}(B_{1}(0))$ non-negative, such that $ \int_{B_{1}(0)}m_{0} \, dx=1$ and $\int_{B_{1}(0)}\frac{\abs{\nabla m_{0}}^{\gamma'}}{m_{0}^{\gamma'-1}} \, dx$ is finite. Now pick $x_{0}\in \Omega$ and define
\begin{equation*}
m_{\lambda}(x)=\lambda^{N}m_{0}(\lambda (x-x_{0})) \qquad w_{\lambda}(x)=\nabla m_{\lambda}(x)=\lambda^{N+1} \nabla m_{0}(\lambda (x-x_{0})).
\end{equation*}
We notice that for $\lambda>\frac{1}{\dist (x_{0},\partial\Omega)}$, then $\supp (m_{\lambda})\subset B_{\frac{1}{\lambda}}(x_{0})$ and $(m_{\lambda},w_{\lambda})\in \mathcal{K}$. Moreover $\norm{m_{\lambda}}_{q}^{q}=\lambda^{N(q-1)}\norm{m_{0}}_{q}^{q}$.
Now we have, using our assumption and \eqref{e:mL:prop},
\begin{multline*}
    \mathcal{E}(m_{\lambda},w_{\lambda})=\nrg{m_{\lambda}}{w_{\lambda}}\leq \\
    \cl^{-1}\int_{\Omega}\frac{\abs{w_{\lambda}}^{\gamma'}}{m_{\lambda}^{\gamma'-1}}\, dx-\frac{\cf}{q}\norm{m_{\lambda}}_{q}^{q}-C=\lambda^{\gamma'}\cl^{-1}\int_{B_{1}(0)}\frac{\abs{\nabla m_{0}}^{\gamma'}}{m_{0}^{\gamma'-1}}\, dt-\lambda^{N(q-1)}\frac{\cf}{q}\norm{m_{0}}_{q}^{q}-C
\end{multline*}
and we see that under our assumptions on $q$ for $\lambda\rightarrow +\infty$, the right-hand side goes to $-\infty$.
By the above computations, we can also see that $\mathcal{E}$ may be unbounded in the case $q=1+\frac{\gamma'}{N}$, provided that $\cf$ is large enough. 
\end{remark}
Let us go back to the existence of a local minimum in the case $q>\Bar{q}$. To do so, we need to impose some restriction on the coefficient $\cf$, so that an $\alpha$ such that $c_{\alpha}<\hat{c}_{\alpha}$ can be found. 
\begin{theorem}
\label{th:minim}
Let $q>\Bar{q}$,
\begin{equation}
\label{e:alpha}
  \Bar{\alpha}=\left(\frac{\cl}{\cf\cpr}\right)^{q'}  ,
\end{equation}
and suppose that
\begin{equation}
 \label{e:cf}
 \cf<\min\left\{\frac{\cl}{\cpr},\left(K'q'\right)^{1-q}\left(\frac{\cl}{\cpr}\right)^{q}\right\},
\end{equation}
where 
\begin{equation}
K':=\cl+\cl^{-1}+2K_{f}+\frac{\cl}{q\cpr}.
\end{equation}
Then, there exists $(m_{\ep},w_{\ep})\in \mathcal{K}\cap B_{\Bar{\alpha}}$
such that $(m_{\ep},w_{\ep})$ is a local minimum for $\mathcal{E}_{\ep}$. Moreover it holds
\begin{equation}
    \label{e:minim}
    \mathcal{E}_{\ep}(m_{\ep},w_{\ep})=c_{\Bar{\alpha}-\delta}=c_{\Bar{\alpha}}
\end{equation}
for all $\delta$ less than some $\Bar{\delta} >0$ small enough. 
\end{theorem}
\begin{proof}
Let us start the by the simple observation that if we find $\alpha_{2}>\alpha_{1} $
such that $\hat{c}_{\alpha_{2}}>\hat{c}_{\alpha_{1}}$, then we have:
\begin{equation*}
    c_{\alpha_{2}}=\min_{\alpha}\{ \hat{c}_{\alpha} : 0 \leq \alpha\leq \alpha_{2} \}\leq \hat{c}_{\alpha_{1}}<\hat{c}_{\alpha_{2}}
\end{equation*}
and so we can conclude the existence of an interior minimum in $B_{\alpha_{2}}$.
We choose $\alpha_{1}=1$ and $\alpha_{2}=\Bar{\alpha}$.
Since $\cf < {\cl} / {\cpr}$, we are sure that $\bar \alpha > 1$.\\
We now consider estimates for $\hat{c}_{\alpha_{1}}$.
We begin by noticing that by H\"older inequality we have
\begin{equation*}
    1=\norm{m}_{1}\leq \norm{m}_{q}\abs{\Omega}^{\frac{1}{q'}}=1
\end{equation*}
which implies $m=1$ a.e.\, and therefore $(m,w)\equiv(1,0)\in \mathcal{K}\cap U_{\alpha_{1}}$. Thus, using \eqref{e:mL:prop} and \eqref{e:Feps:ass} we find
\begin{equation}
\label{e:alphaupper}
    \hat{c}_{\alpha_{1}}\leq \epsen{1}{0}\leq \cl^{-1}+\frac{\cf}{q}+K_{f}<\cl^{-1}+\frac{\cl}{q\cpr}+K_{f}.
\end{equation}
thanks to our assumptions on $\cf$.
Next, we rewrite \eqref{e:energybndal} with $\alpha=\Bar{\alpha}$ as
\begin{equation}
\label{e:alphalower}
    \hat{c}_{\Bar{\alpha}}\geq \frac{\cl}{\cpr}\Bar{\alpha}^{\frac{1}{q}}-\cl-\frac{\cf}{q}\Bar{\alpha}-K_{f}.
\end{equation}

To conclude that $\hat{c}_{\bar \alpha}>\hat{c}_{1}$, we need to check that 
\begin{equation}
\label{e:alphaconst}
    \phi(\alpha):=\frac{\cl}{\cpr}\alpha^{\frac{1}{q}}-\frac{\cf}{q}\alpha > \cl+\cl^{-1}+2K_{f}+\frac{\cl}{q\cpr}=K',
\end{equation}
and the previous inequality holds again by the assumptions on $C_f$ (notice that $\Bar{\alpha}$ maximizes $\phi$). Moreover since $\phi$ is continuous, we get that some $\delta$ small enough exists so that also $\hat{c}_{\bar{\alpha}-\delta}>\hat{c}_{\alpha_{1}}$.
\end{proof}
In the previous construction, the assumptions on $\cf$ are chosen to have the largest possible $\alpha$ such that the energy admits an interior minimizer in $B_{\bar \alpha}$, and $\Bar{\alpha}$ depends on the value of $\cf$. To treat the case $q = \qc$, we will need to find a minimizer in $B_{\Bar{\alpha}}$ independent of $\cf$. This is possible under different assumptions on $\cf$. 
\begin{theorem}
\label{th:minim2}
Let
\begin{equation}
    \hat{\alpha}=\left(\frac{\cpr}{\cl}K''+1\right)^q
\end{equation}
where
\begin{equation*}
    K''=\cl+\cl^{-1}+2K_{f}
\end{equation*}
and suppose 
\begin{equation}
\label{e:cfq*}
    \cf<\frac{q\cl}{\cpr(\hat{\alpha}+1)}
\end{equation}
then, the results of Theorem \ref{th:minim} hold. 
\end{theorem}
\begin{proof}
The proof is analogous to the one of Theorem \ref{th:minim}. The only differences are the choice of  $\alpha$ when enforcing \eqref{e:alphaconst} ($\hat{\alpha}$ no longer maximizes $\phi(\alpha)$) and the last inequality in \eqref{e:alphaupper}, which is skipped.
\end{proof}

\subsection{Convex duality}
We now employ some convex duality arguments to obtain, from the (local) minimizer constructed in the previous section, a solution to the MFG system \eqref{e:sysreg}. We follow the usual route (see e.g. \cite{cesaroni_cirant_2019} and references therein), which requires first to linearize the functional (which is not convex by the presence of the possibly nonconvex $F_\eps$) around the minimizer that we found.
Given $(m_{\ep},w_{\ep})$, which is a global minimizer of $\mathcal{E}_{\ep}$ on $\mathcal{K}$, or a local minimizer on $\mathcal{K}\cap B_{\Bar{\alpha}}$ when $q > \bar q$, let us introduce the following linearized functional:
\begin{equation}
\label{e:jeps}
J_{\ep}(m,w)=\int_{\Omega}\ml{m}{w}+f_{\ep}[m_{\ep}](x)m \, dx.
\end{equation}
We notice that this functional is convex. We now prove that this functional admits the same minimizer as $\mathcal{E}_{\ep}$.
\begin{proposition}
\label{prop:samemin}
Let $(m_{\ep},w_{\ep})$ be a global minimizer of $\mathcal{E}_{\ep}$ on $\mathcal{K}$ or a local minimizer on $\mathcal{K}\cap B_{\Bar{\alpha}}$ as constructed above. Then
\begin{equation}
    \label{e:samemin}
     \min_{(m,w)\in \mathcal{K}}J_{\ep}(m,w)=J_{\ep}(m_{\ep},w_{\ep}).
\end{equation}
\end{proposition}
\begin{proof}
Let $(m,w)\in \mathcal{K}$ and consider for $0<\lambda<1$
\begin{equation*}
    m_{\lambda}=\lambda m +(1-\lambda)m_{\ep}
\end{equation*}
If $(m_{\ep},w_{\ep})$ is a local minimum, since by \eqref{e:minim} $m_{\ep}\in B_{\Bar{\alpha}-\delta}$ for some positive $\delta$, we can conclude that for $\lambda$ small enough, $m_{\lambda}\in B_{\Bar{\alpha}}$. If $(m_{\ep},w_{\ep})$ is a global minimum this argument holds for all $\lambda$. 
Hence, by minimality and convexity 
\begin{multline*}
    F_{\ep}[m_{\ep}]-F_{\ep}[m_{\lambda}]\leq
    \int_{\Omega}\ml{m_{\lambda}}{w_{\lambda}}\, dx-\int_{\Omega}\ml{m_{\ep}}{w_{\ep}}\, dx \\ \leq \lambda \int_{\Omega}\ml{m}{w}\, dx+(1-\lambda)\int_{\Omega}\ml{m_{\ep}}{w_{\ep}}\, dx-\int_{\Omega}\ml{m_{\ep}}{w_{\ep}}\, dx  \\ =
    \lambda \left(\int_{\Omega}\ml{m}{w}\, dx - \int_{\Omega}\ml{m_{\ep}}{w_{\ep}}\, dx \right).
\end{multline*}
Now, by \eqref{e:l2der}, we have
\begin{multline*}
 F_{\ep}[m_{\ep}]-F_{\ep}[m_{\lambda}]=\int_{0}^{1}\int_{\Omega}f_{\ep}[(1-t)m_{\lambda}+tm_{\ep}](x)(m_{\ep}-m_{\lambda})(x)\, dx\\=
 -\lambda \int_{0}^{1}\int_{\Omega}f_{\ep}[m_{\ep}+\lambda (1-t)(m-m_{\ep})](x)(m-m_{\ep})(x)\, dx.
\end{multline*}
Combining the two expressions, we can use Lipschitz estimates for $f_{\ep}[\cdot](x)$ near $m_{\ep}$ and send $\lambda$ to $0$ to conclude that
\begin{equation*}
\int_{\Omega}\ml{m}{w}\, dx-\int_{\Omega}\ml{m_{\ep}}{w_{\ep}}\, dx
\geq-\int_{\Omega}f_{\ep}[m_{\ep}](x)(m_{\ep}-m)(x)\, dx,
\end{equation*}
which is equivalent to the minimality of $J_{\ep}$ (globally on $\mathcal K$).
\end{proof}
Now that we have global minimizer of a convex functional, we can construct a solution of \eqref{e:sysreg}.
\begin{theorem}
\label{th:duality}
Let $(m_{\ep},w_{\ep})$ be a minimizer of $J_{\ep}$ constructed above. Then 
$m_{\ep}\in W^{1,p}(\Omega)$ for all $p>1$ and there exists $\lambda_{\ep}\in \mathbb{R}$ and $u_{\ep}\in C^{2}(\overline{\Omega})$ such that $(u_{\ep},\lambda_{\ep},m_{\ep})$ is a solution to \eqref{e:sysreg}. Moreover,
\begin{equation}
\label{e:wcara}
    w_{\ep}=-m_{\ep}\nabla H(\nabla(u_{\ep})),
\end{equation}
and there exists $C>0$ independent of $\ep$ such that 
\begin{equation}
\label{e:mbounds}
    \norm{m_{\ep}}_{q}\leq C, \qquad \norm{m_{\ep}}_{1,r}\leq C,
\end{equation}
and
\begin{equation}
\label{e:lambdabnd}
     \abs{\lambda_{\ep}}\leq C.
\end{equation}
\end{theorem}
\begin{proof}
The proof is like \cite[Th.\ 4]{cesaroni_cirant_2019} with minor modifications.
\end{proof}
\subsection{Passage to the limit}
We now wish to let $\ep\rightarrow 0$, and to do so we need some a priori estimate. We distinguish two cases: if $q<\qc$, we can use a blow-up argument to deduce an a priori $L^{\infty}$ bound on $m_{\ep}$ (that by a bootstrap procedure yields further estimates on $u,m$ and their derivatives). If $q=\qc$, the argument fails and we need to require some extra smallness on $\cf$ to obtain such bound. 

The blow up argument follows the lines of \cite{cesaroni_cirant_2019}, but we need an extra care for the presence of the Neumann boundary conditions.
\begin{proposition}
\label{prop:apriorilinf}
Let $(u_{\ep},\lambda_{\ep},m_{\ep})$ be a solution to \eqref{e:sysreg} constructed above, and suppose that $q<\qc$.
Then there exist $C>0$ independent of $\ep$ such that
\begin{equation}
\label{e:apriorilinf}    
\norm{m_{\ep}}_{\infty}\leq C.
\end{equation}
\end{proposition}
\begin{proof}
Suppose by contradiction that
\begin{equation*}
    M_{\ep}=\max_{\Omega}m_{\ep}=m_{\ep}(x_{\ep})\rightarrow + \infty.
\end{equation*}
Define 
\begin{equation*}
    \mu_{\ep}:=M_{\ep}^{-\beta} \qquad \beta:=(q-1)\frac{\gamma-1}{\gamma}.
\end{equation*}
We have that $\mu_{\ep}\rightarrow 0$.
Define the following rescaling
\begin{equation*}
    \begin{cases}
    v_{\ep}(x)&=\mu_{\ep}^{\frac{2-\gamma}{\gamma-1}}u_{\ep}(\mu_{\ep}x+x_{\ep})\\
    n_{\ep}(x)&=M_{\ep}^{-1}m_{\ep}(\mu_{\ep}x+x_{\ep}).
    \end{cases}
\end{equation*}
We notice that $n_{\ep}(0)=1$ and that $0\leq n_{\ep}(x)\leq 1$.
Define also
\begin{equation}
    H_{\ep}(q)=\mu_{\ep}^{\frac{\gamma}{\gamma-1}}H(\mu_{\ep}^{\frac{1}{1-\gamma}}q)\qquad \nabla H_{\ep}(q)=\mu_{\ep}\nabla H(\mu_{\ep}^{\frac{1}{1-\gamma}}q)
\end{equation}
By \eqref{e:ass:H} we get
\begin{equation}
\label{e:Hrescaled}
    C^{-1}_{H}\abs{p}^\gamma-K_{H}\leq H_{\ep}(p)\leq C_{H}\abs{p}^\gamma+K_H \qquad 
     \abs{\nabla H_{\ep}(p)}\leq C_{H}\abs{p}^{\gamma -1} + K_H.
\end{equation}
Then, define
\begin{equation}
    \Tilde{f_{\ep}}(x):=\mu_{\ep}^{\frac{\gamma}{\gamma-1}}f_{\ep}[m_{\ep}](x_{\ep}+\mu_{\ep}x).
\end{equation}
Since $m_{\ep}(x)\leq M_{\ep}$, we can use \eqref{e:ass:f} to get that
\begin{multline*}
    \norm{\Tilde{f_{\ep}}[m_{\ep}]}_{\infty}\leq \mu_{\ep}^{\frac{\gamma}{\gamma-1}}\norm{f\left(\cdot,m*\chi_{\ep}(\cdot)\right)*\chi_{\ep}}_{\infty}\leq \mu_{\ep}^{\frac{\gamma}{\gamma-1}}\norm{f\left(\cdot,m_{\ep}*\chi_{\ep}\right)}_{\infty} \\ \leq \mu_{\ep}^{\frac{\gamma}{\gamma-1}} C\left(\norm{m_{\ep}*\chi_{\ep}}^{q-1}_{\infty}+1\right)\leq \mu_{\ep}^{\frac{\gamma}{\gamma-1}}\left( C\norm{m_{\ep}}^{q-1}_{\infty}+C\right)  \\ \leq C+CM_{\ep}^{q-1-\beta \frac{\gamma}{\gamma-1}}\leq C
\end{multline*}
by our definitions of $\mu_{\ep}$ and $\beta$. 
Lastly, define 
\begin{equation}
    \Tilde{\lambda}_{\ep}=\mu_{\ep}^{\frac{\gamma}{\gamma-1}}\lambda_{\ep}.
\end{equation}
Clearly by \eqref{e:lambdabnd}, $\abs{\Tilde{\lambda}_{\ep}}\leq C$.
Now after some computations we have
\begin{equation*}
    \begin{cases}
    \Delta v_{\ep}(x)=\mu_{\ep}^{\frac{\gamma}{\gamma-1}}\Delta u_{\ep}(x\mu_{\ep}+x_{\ep})\\
    H_{\ep}(\nabla v_{\ep}(x))=\mu_{\ep}^{\frac{\gamma}{\gamma-1}}H(\nabla 
    u_{\ep}(x\mu_{\ep}+x_{\ep}))\\
    \Delta n_{\ep}(x)=\mu_{\ep}^{\frac{1}{\beta}+2}\Delta m_{\ep}(x\mu_{\ep}+x_{\ep})\\
    \nabla H_{\ep}(\nabla v_{\ep}(x))=\mu_{\ep}\nabla H(\nabla u_{\ep}(x\mu_{\ep}+x_{\ep}))\\
    \dive(n_{\ep}\nabla H_{\ep}(\nabla v_{\ep}(x)))=\mu_{\ep}^{\frac{1}{\beta}+2}\dive(m_{\ep}(x\mu_{\ep}+x_{\ep})\nabla H(\nabla u_{\ep}(x\mu_{\ep}+x_{\ep})))\\
     \frac{\partial v_{\ep}(x)}{\partial n }=\mu_{\ep}^{\frac{1}{\gamma-1}}\frac{\partial u_{\ep}(\mu_{\ep}x+x_{\ep})}{\partial n }(\mu_{\ep}x+x_{\ep})\\
    \frac{\partial n_{\ep} }{\partial n}=\mu_{\ep}^{\frac{1}{\beta}+1}\frac{\partial m_{\ep}(\mu_{\ep}x+x_{\ep}) }{\partial n}.
    \end{cases}
\end{equation*}
from which we deduce that $(v_{\ep},\Tilde{\lambda}_{\ep},n_{\ep})$ is a solution of 
\begin{equation}
\label{e:sysrescaled}
\begin{cases}
-\Delta v_{\ep} +H_{\ep}(\nabla v_{\ep} )+ \Tilde{\lambda}_{\ep} = \Tilde{f_{\ep}}(x) \qquad & \text{on} \ \Omega_{\ep} \\
-\Delta n_{\ep} - \dive({n_{\ep}\nabla H_{\ep} (\nabla v_{\ep} )})=0 \qquad  & \text{on} \ \Omega_{\ep} \\
\frac{\partial v_{\ep}}{\partial n }=0 \qquad & \text{on} \  \partial \Omega_{\ep} \\
\frac{\partial n_{\ep}}{\partial n } + n_{\ep} \nabla H_{\ep}(\nabla v_{\ep} )\cdot n =0 \qquad & \text{on} \ \partial \Omega_{\ep}\\
\int_{\Omega_{\ep}}n_{\ep}=M_{\ep}^{-1-\beta}, \qquad \int_{\Omega}v_{\ep}=0.
\end{cases}
\end{equation}
where $\Omega_{\ep}=\{x:\mu_{\ep}x+x_{\ep}\in \Omega\}$.
We now have to distinguish two cases. 
Suppose first that
\begin{equation*}
    \lim_{\ep\rightarrow0}\frac{d\left(x_{\ep},\partial\Omega\right)}{\mu_{\ep}}=+\infty.
\end{equation*}
From that we have that $\Omega_{\ep}\uparrow\rn$, hence for $\ep$ small enough we have that $\Omega_{\ep}\supset B_{4R}(0)$ for an $R>0$ independent on $\eps$. 
We know that $\Tilde{\lambda}_{\ep}$ and $\Tilde{f}_{\ep}(x)$ are uniformly bounded, thus using Theorem \ref{th:hjexistence} (and the remark below) we can conclude that there exists $C$ independent of $\ep$ such that $\norm{\nabla v_{\ep}}_{\infty}\leq C$ on $B_{2R}$. Now, using \eqref{e:Hrescaled} we can deduce that $\norm{n_\ep \nabla H_{\ep}(\nabla v_{\ep})}_{\infty}\leq C$. Hence, by Proposition \ref{prop:lpest}, for $\ep$ small enough $n_{\ep}$ is equibounded in $W^{1,p}(B_{R}(0))$ for all $p>1$, and by Sobolev Embedding also in $C^{\theta}(\overline{B}_{R})$.
We know that $n_{\ep}(0)=1$, therefore using equiboundness in $C^{\theta}(\overline{B}_{R})$ we can deduce that there exists $\delta>0$ and $r<R$ such that $\int_{B_{r}(0)}n_{\ep}^{q}(x)\, dx>\delta>0$. Thus we have that
\begin{equation}
\label{e:blowupcontr}
    0<\delta \leq \int_{B_{r}(0)}n_{\ep}^{q}(x)\, dx\leq \norm{n_{\ep}}_{q}^{q}=M_{\ep}^{-q}\mu_{\ep}^{-N}\norm{m_{\ep}}_{q}^{q}=M_{\ep}^{-q+\beta N}\norm{m_{\ep}}_{q}^{q}.
\end{equation}
Since $q<\qc$ then $-q+\beta N<0$, and using \eqref{e:mbounds} we have that
\begin{equation*}
    0< \delta \leq M_{\ep}^{-q+\beta N}\norm{m_{\ep}}_{q}^{q}\leq C M_{\ep}^{-q+\beta N}\rightarrow 0 
\end{equation*}
which is a contradiction.

Suppose now that
\begin{equation*}
    \lim_{\ep\rightarrow0}\frac{d\left(x_{\ep},\partial\Omega\right)}{\mu_{\ep}}\leq C.
\end{equation*}
Up to subsequences we can suppose that $x_{\ep}\rightarrow \bar x \in \partial\Omega$ as $\ep\rightarrow 0$. Moreover, up to an affine transformation we can assume $\bar x=0\in\partial\Omega$ and $n(0)=-e_{N}$. Define $x'=(x_{1},\dots,x_{N-1})$. By the smoothness of $\Omega$ there exists $U\subset \rn$, $\Gamma\subset \mathbb{R}^{N-1}$ and $\phi(x') \in C^{2,\alpha}(\Gamma) $ such that 
\begin{align*}
   &\phi(0)=0, \qquad \nabla \phi (0)=0,\\
   &\partial\Omega\cap U =\{ (x',x_{N}) : x_{N}=\phi(x')\},\\
   &\Omega\cap U=\{ (x',x_{N}) : x_{N}>\phi(x')\}.
\end{align*}
Let us now define a diffeomorfism $\Psi:\rn\rightarrow\rn$ that "straightens" the boundary. We set 
\begin{equation}
   y_{i}= (\Psi(x))_{i}:=
    \begin{cases}
    x_{i}-x_{N}\frac{\partial\phi}{\partial x_{i}}(x') \qquad &\text{ for } 1\leq i\leq N-1\\
    x_{N}-\phi(x') \qquad &\text{ for } i=N.
    \end{cases}
\end{equation}
We can see that $\Psi$ is invertible in a neighborhood of $0$. 
We now extend with an even reflection $v_{\ep}$ and $m_{\ep}$. We set 
\begin{align}
    & w_{\ep}(y)= v_{\ep}\left(\frac{\Psi^{-1}(y',\abs{y_{N}})-x_{\ep}}{\mu_{\ep}}\right)\\
    & \rho_{\ep}(y)= n_{\ep}\left(\frac{\Psi^{-1}(y',\abs{y_{N}})-x_{\ep}}{\mu_{\ep}}\right).
\end{align}
Due to the homogeneous Neumann boundary conditions, with some calculations it is possible to show that $\frac{\partial w_{\ep}}{\partial y_{N}}\vert_{\{y_{N}=0\}}=0$. Moreover, one can derive that $w_{\ep},\rho_{\ep}$ satisfy an equation similar to \eqref{e:sysrescaled} in a fixed neighborhood of the boundary point $p$ independent of $\eps$ (with coefficients that converge to the identity as $\eps \to 0$). From this, we can repeat the above argument and reach a contradiction. 
\end{proof} 
We can see in this proof the criticality of the case $q=\qc$. Looking at \eqref{e:blowupcontr} and the lines below, in the case $q=\qc$ it is not possible to reach a contradiction, since $M_{\ep}^{-q+\beta N}$ does not vanish. 
To tackle this problem, the idea is to obtain additional regularity using finer estimates for both the Fokker-Planck and the Hamilton-Jacobi equation. Once we find uniform bounds for $m_{\eps}$ in some $L^{p}$ with $p>\qc$, the above arguments can be used to conclude again that we have an uniform $L^{\infty}$ bound. This procedure requires additional assumptions on $\cf$. Moreover we need to have a value for $\Bar{\alpha}$ which is independent from $\cf$, as we have constructed in Theorem \ref{th:minim2}. 

Below, $C_S$ is the Sobolev Embedding constant for $W^{1,p}(\Omega)$ into $L^{p^*}(\Omega)$, where $p^*$ is chosen so that
\[
p^*<N \quad \text{and} \quad p^* > 1 + \frac{\gamma'}{N-\gamma'}.
\]

\begin{proposition}
\label{prop:estcritical}
Let $(u_{\ep},\lambda_{\ep},m_{\ep})$ be a solution to \eqref{e:sysreg} and suppose that $q=\qc$.
Then 
\begin{equation}
\label{e:aprioriln}    
\norm{\nabla u_{\ep}}_{N(\gamma-1)}\leq \frac{1}{\left(4C_{E}C_{S}C_{H}\right)^{\frac{1}{\gamma-1}}},
\end{equation}
where $C_{E}$ is defined in  Proposition \ref{prop:lpest}, provided that $C_f$ and $K_f$ are small enough (that is, smaller than some positive constant depending on $\Omega, N, \gamma, q$).
\end{proposition}
\begin{proof}
We use Theorem \ref{th:hjapriori2}, choosing $\delta$ small enough so that $M(\delta)\le\left(4C_{E}C_{S}C_{H}\right)^{\frac{-1}{\gamma-1}}$. Let us compute the norm of $f_{\ep}$ in $\leb{\frac{N}{\gamma'}}$. Using \eqref{e:ass:f}, convolution properties and the definition of $\qc$ we have
\begin{equation*}
    \norm{f_{\ep}}_{\frac{N}{\gamma'}}=\norm{f(x, m*\chi_{\ep}(x))}_{\frac{N}{\gamma'}}\leq \cf \norm{m_{\ep}^{q-1}}_{\frac{N}{\gamma'}}+K_{f}=\cf(\norm{m_{\ep}}^{q}_{q})^\frac{1}{q'}+K_{f}.
\end{equation*}
Now using Theorem \ref{th:minim2}, we know that 
\begin{equation*}
    \norm{m}_{q}^{q}\leq \hat{\alpha}= \left(\frac{\cpr}{\cl}K''+1\right)^q.
\end{equation*}
Hence if
\begin{equation*}
    \cf \hat{\alpha}^{\frac{1}{q'}}\ + K_f \leq \delta
\end{equation*}
we have that $\norm{f_{\ep}}_{\frac{N}{\gamma'}}\leq \delta$. Thus, we can apply Theorem \ref{th:hjapriori2}  to conclude.
\end{proof}
\begin{corollary}
\label{cor:estcritical}
Let $(u_{\ep},\lambda_{\ep},m_{\ep})$ be a solution to \eqref{e:sysreg} and suppose that $q=\qc$. Under the assumptions of the previous proposition, supposing in addition that
\begin{equation}\label{e:khs}
K_H \le \frac{1}{4C_{E}C_{S}},
\end{equation}
 then there exist $C>0$ independent of $\ep$ such that
\[
\norm{m_{\ep}}_{\infty}\leq C.
\]
\end{corollary}
\begin{proof} Using Proposition \ref{prop:estcritical}
we have that $\norm{\nabla u_{\ep}}_{N(\gamma-1)}\leq \frac{1}{\left(2 C_{E}C_{S}C_{H}\right)^{\frac{1}{\gamma-1}}}$. Hence $$\norm{\nabla H(\nabla u_{\ep})}_{N}\leq \frac{1}{4C_{E}C_{S}}+K_H$$ and we can use Proposition \ref{prop:fpestcrit} to conclude that $m_{\ep}$ are uniformly bounded in $W^{1,p}(\Omega)$ for the chosen above; in this way,  by Sobolev Embeddings $m_\ep$ are uniformly bounded in $\leb{q}$ for some $q>\qc$. Once we have this bound, we can proceed as in Proposition \ref{prop:apriorilinf} to conclude that $m_{\ep}$ is bounded in $L^{\infty}$.
\end{proof}

Now everything is ready to prove the main result. Thanks to the uniform bounds, we use a bootstrap procedure to obtain the regularity which is necessary to pass to the limit into the equations. 
\begin{proof} [Proof of Theorem \ref{th:existence}] 
We set ourselves in the assumptions of the previous propositions and theorems, in particular we require $C_f$, $K_f$ and $K_H$ to be possibly small enough (see in particular Lemma \ref{le:minim}, equations \eqref{e:cf}, \eqref{e:cfq*}, Proposition \ref{prop:estcritical}). 

Let $(u_{\ep},\lambda_{\ep},m_{\ep})$ be the sequence of solutions constructed in Theorem \ref{th:duality}. We first show that we can pass to the limit as $\eps\to0$ and obtain a solution 
of \eqref{e:sys}; in this we need to treat separately the Sobolev subcritical and critical cases. Next, we will prove that the variational characterization passes to the limit as well.

Let us first suppose $q<\qc$. By \eqref{e:lambdabnd} we have that, up to subsequences, $\lambda_{\ep}\rightarrow \lambda$. Moreover, by Proposition \ref{prop:apriorilinf} we have that $\norm{m_{\ep}}_{\infty}\leq C$. Hence by properties of mollifiers and \eqref{e:ass:f} we get
\begin{multline*}
    \norm{f_{\ep}[m_{\ep}]}_{\infty}=\norm{f\left(\cdot,m*\chi_{\ep}(\cdot)\right)*\chi_{\ep}}_{\infty}\leq \norm{f\left(\cdot,m_{\ep}*\chi_{\ep}\right)}_{\infty} \\\leq C\norm{m_{\ep}*\chi_{\ep}}^{q-1}_{\infty}+C\leq C\norm{m_{\ep}}^{q-1}_{\infty}+C\leq C
\end{multline*}
Thus, we can conclude by Theorem \ref{th:hjexistence}, that $\norm{\nabla u_{\ep}}_{\infty}\leq K$ for some $K>0$ independent of $\ep$. Using  the estimates on $H$ and elliptic regularity in the Hamilton-Jacobi equation we can conclude that $\norm{u_{\ep}}_{1,p}\leq C$ for all $p>1$. Moreover, by Sobolev embeddings,
$u_{\ep}$ is equibounded in $C^{1,\alpha}(\overline{\Omega})$ for all $\alpha\in (0,1)$. Since $\nabla H(\nabla u_{\ep})\leq C$, we now use Theorem \ref{th:fp}
to conclude that  $\norm{m_{\ep}}_{1,p}\leq C$ for all $p>1$, and therefore by Sobolev embedding $m_{\ep}$ is equibounded in $C^{\theta}(\overline{\Omega})$ for all $\theta\in (0,1)$. Hence, we get that, up to subsequences, $m_{\ep}\rightharpoonup m$ in $W^{1,p}$ for all $p>1$ and $m_{\ep}\rightarrow m$ uniformly. We can then go back to the Hamilton-Jacobi equation and, with a similar reasoning, conclude that $f_{\ep}[m_{\ep}](x)$ is equibounded in $C^{\theta}(\Omega)$ for all $\theta\in(0,1)$. Hence $u_{\ep}$ is equibounded in $C^{2,\theta}(\overline{\Omega})$ for all $\theta\in(0,1)$. Finally, we can conclude that up to subsequences $u_{\ep}\rightarrow u$ in $C^{2}(\overline{\Omega})$. Now the convergences are strong enough to pass to the limit in the equations, so we can conclude that $(u,\lambda,m)$ is a solution of \eqref{e:sys}, with the positivity of $m$ coming from Theorem \ref{th:fp} and pointwise convergence. 

When $q=\qc$, we argue in the very same way, starting from Corollary \ref{cor:estcritical}.

Finally, we are left to prove that the solutions we found are (local) minimizers of $\Ecal$. We will use the fundamental theorem of $\Gamma-$convergence (see e.g.\ \cite{braides_2002}). Notice that, since we know a priori that the sequence of minima converges, we do not need to prove an equicoercivity result. Let us show that $\mathcal{E}_{\ep}$ $\Gamma-$converges to $\mathcal{E}$ on the space $X=L^{q}(\Omega)\cap {W}^{1,r}(\Omega)\times\leb{1}$. Suppose that $(m_{\ep},w_{\ep})\rightarrow (m,w)$ in $X$. By properties of mollifiers and continuity of the convolution we have $m_{\eps}*\chi_{\eps}\rightarrow m$ in $\leb{q}$. We already remarked the semicontinuity of the Lagrangian term in $\mathcal{E}_{\eps}$, moreover by \eqref{e:ass:f} we have strong $\leb{q}$ continuity of $m \rightarrow \int_{\Omega}\int_{0}^{m}f(x,n)\,dn\, dx$. Thus
\begin{multline*}
        \liminf_{\ep} \mathcal{E}_{\ep}(m_{\ep},w_{\ep})= \liminf_{\ep}\int_{\Omega}\ml{m_{\ep}}{w_{\ep}}\,dx+\lim_{\eps}F_{\eps}[m_{\eps}]\\ \geq \int_{\Omega}\ml{m}{w}\,dx+ \int_{\Omega}F(x,m)\,dx=\mathcal{E}(m,w).
\end{multline*}
As for the recovery sequence, it suffices to choose $(m_{\ep},w_{\ep})=(m,w)$ for all $\ep>0$ and we clearly have that $\mathcal{E}_{\ep}(m,w)\rightarrow \mathcal{E}(m,w)$ by the properties of mollifiers and again by the strong $\leb{q}$ continuity. 
To finish, we know by \eqref{th:duality} that a minimum $(m_{\ep},w_{\ep})$ of $\mathcal{E}_{\ep}$ yields a solution $(u_{\ep},\lambda_{\ep},m_{\ep})$ of \eqref{e:sysreg} and that the relation $w_{\ep}=-m_{\ep}\nabla H(\nabla(u_{\ep}))$ holds. Since we know that  $(u_{\ep},\lambda_{\ep},m_{\ep})$ converges in $C^{2}(\overline{\Omega})\times \R\times W^{1,p}(\Omega)$ for all $p$ to a solution $(u,\lambda,m)$ of the original problem, we get that $(m_{\ep},w_{\ep})$ converges in $X$ to $(u,-m\nabla H(\nabla u ))$.
Hence, we can conclude that the solution $(u,\lambda,m)$ is such that $(m,-m\nabla H(\nabla u))$ is a minimum of $\mathcal{E}$ (possibly restricted to $B_{\hat \alpha}$ when $q > \bar q$).
\end{proof}

\bigskip

\textbf{Data Availability.} Data sharing not applicable to this article as no datasets were generated or analyzed during the current study.

\bigskip

\begin{wrapfigure}{r}{0.15\linewidth}
\includegraphics[width=\linewidth]{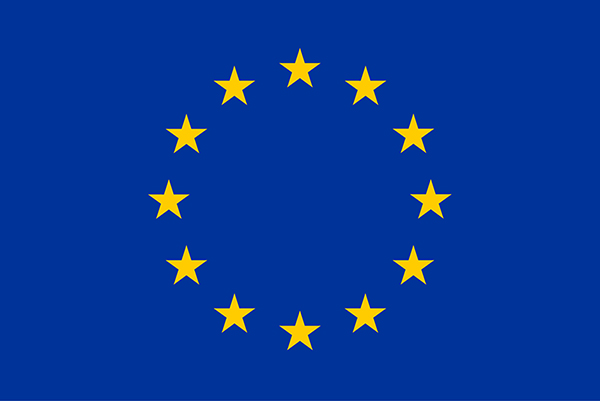}
\end{wrapfigure}
\textbf{Acknowledgements.} M.C. is partially supported by the King Abdullah University of Science and Technology (KAUST) project CRG2021-4674 ``Mean-Field Games: models, theory, and computational aspects''; A.C. is partially supported by the European Union's Horizon 2020 research and innovation programme under the Marie Sklodowska-Curie grant agreement No 945332; G.V. is partially supported by the project Vain-Hopes within the program  VALERE-Universit\`a degli Studi della Campania ``Luigi Vanvitelli'', by the Portuguese 
government through FCT/Portugal under the project PTDC/MAT-PUR/1788/2020; M.C. and G.V. are members of the Gruppo Nazionale per l'Analisi Matematica, la Probabilit\`a e le loro Applicazioni (GNAMPA) of the Istituto Nazionale di Alta Matematica (INdAM).

\bibliographystyle{abbrv}
\bibliography{Bibliography}{}

\begin{thebibliography}{10}

\bibitem{agmon_1959}
S.~Agmon.
\newblock The {$L_{p}$} approach to the {D}irichlet problem. {I}. {R}egularity
  theorems.
\newblock {\em Ann. Scuola Norm. Sup. Pisa Cl. Sci. (3)}, 13:405--448, 1959.

\bibitem{MR2816471}
A.~Ambrosetti and D.~Arcoya.
\newblock {\em An introduction to nonlinear functional analysis and elliptic
  problems}, volume~82 of {\em Progress in Nonlinear Differential Equations and
  their Applications}.
\newblock Birkh\"{a}user Boston, Ltd., Boston, MA, 2011.

\bibitem{bardi_feleqi}
M.~Bardi and E.~Feleqi.
\newblock Nonlinear elliptic systems and mean field games.
\newblock {\em NoDEA Nonnlinear Differ. Equ. Appl.}, 23:23--44, 2016.

\bibitem{bensoussan_1988}
A.~Bensoussan.
\newblock {\em Perturbation methods in optimal control}.
\newblock Wiley/Gauthier-Villars Series in Modern Applied Mathematics. John
  Wiley \& Sons, Ltd., Chichester; Gauthier-Villars, Montrouge, 1988.
\newblock Translated from the French by C. Tomson.

\bibitem{bernardini_cesaroni}
C.~Bernardini and A.~Cesaroni.
\newblock Ergodic {M}ean-{F}ield {G}ames with aggregation of {C}hoquard-type,
  2022.

\bibitem{braides_2002}
A.~Braides.
\newblock {\em {$\Gamma$}-convergence for beginners}, volume~22 of {\em Oxford
  Lecture Series in Mathematics and its Applications}.
\newblock Oxford University Press, Oxford, 2002.

\bibitem{cesa_cirant_2019}
A.~Cesaroni and M.~Cirant.
\newblock Concentration of ground states in stationary mean-field games
  systems.
\newblock {\em Anal. PDE}, 12(3):737--787, 2019.

\bibitem{cesaroni_cirant_2019}
A.~Cesaroni and M.~Cirant.
\newblock Introduction to variational methods for viscous ergodic mean-field
  games with local coupling.
\newblock In {\em Contemporary research in elliptic {PDE}s and related topics},
  volume~33 of {\em Springer INdAM Ser.}, pages 221--246. Springer, Cham, 2019.

\bibitem{cirant_2014_II}
M.~Cirant.
\newblock On the solvability of some ergodic control problems in {$\Bbb R^d$}.
\newblock {\em SIAM J. Control Optim.}, 52(6):4001--4026, 2014.

\bibitem{cirant_2014}
M.~Cirant.
\newblock Multi-population mean field games systems with {N}eumann boundary
  conditions.
\newblock {\em J. Math. Pures Appl. (9)}, 103(5):1294--1315, 2015.

\bibitem{cirant_2016}
M.~Cirant.
\newblock Stationary focusing mean-field games.
\newblock {\em Comm. Partial Differential Equations}, 41(8):1324--1346, 2016.

\bibitem{cirant_goffi_2021}
M.~Cirant and A.~Goffi.
\newblock On the problem of maximal ${L}^q$-regularity for viscous
  {H}amilton-{J}acobi equations.
\newblock {\em Archive for Rational Mechanics and Analysis}, 240(3):1521--1534,
  mar 2021.

\bibitem{cir_goffi_leonori}
M.~Cirant, A.~Goffi, and T.~Leonori.
\newblock Gradient estimates for quasilinear elliptic {N}eumann problems with
  unbounded first-order terms, 2022.

\bibitem{cirant_porretta}
M.~Cirant and A.~Porretta.
\newblock Long time behaviour and turnpike solutions in mildly non-monotone
  {M}ean {F}ield {G}ames.
\newblock {\em ESAIM: COCV}, 27(86):40, 2021.

\bibitem{cirant_verzini_2017}
M.~Cirant and G.~Verzini.
\newblock Bifurcation and segregation in quadratic two-populations mean field
  games systems.
\newblock {\em ESAIM Control Optim. Calc. Var.}, 23(3):1145--1177, 2017.

\bibitem{cir_verzini_2022}
M.~Cirant and G.~Verzini.
\newblock Local {H}{\"o}lder and maximal regularity of solutions of elliptic
  equations with superquadratic gradient terms.
\newblock {\em Advances in Mathematics}, 409:108700, 2022.

\bibitem{evans_1990}
L.~C. Evans.
\newblock {\em Weak convergence methods for nonlinear partial differential
  equations}, volume~74 of {\em CBMS Regional Conference Series in
  Mathematics}.
\newblock Published for the Conference Board of the Mathematical Sciences,
  Washington, DC; by the American Mathematical Society, Providence, RI, 1990.

\bibitem{GPV}
D.~A. Gomes, E.~A. Pimentel, and V.~Voskanyan.
\newblock {\em Regularity theory for mean-field game systems}.
\newblock SpringerBriefs in Mathematics. Springer, [Cham], 2016.

\bibitem{huang_malhame_caines_2006}
M.~Huang, R.~P. Malham\'{e}, and P.~E. Caines.
\newblock Large population stochastic dynamic games: closed-loop
  {M}c{K}ean-{V}lasov systems and the {N}ash certainty equivalence principle.
\newblock {\em Commun. Inf. Syst.}, 6(3):221--251, 2006.

\bibitem{lasry_lions_2007_II}
J.-M. Lasry and P.-L. Lions.
\newblock Mean field games.
\newblock {\em Jpn. J. Math.}, 2(1):229--260, 2007.

\bibitem{lions_1985}
P.-L. Lions.
\newblock Quelques remarques sur les probl\`emes elliptiques quasilin\'{e}aires
  du second ordre.
\newblock {\em J. Analyse Math.}, 45:234--254, 1985.

\bibitem{noris_tavares_verzini}
B.~Noris, H.~Tavares, and G.~Verzini.
\newblock Normalized solutions for nonlinear {S}chr\"{o}dinger systems on
  bounded domains.
\newblock {\em Nonlinearity}, 32(3):1044--1072, 2019.

\bibitem{MR4191345}
B.~Pellacci, A.~Pistoia, G.~Vaira, and G.~Verzini.
\newblock Normalized concentrating solutions to nonlinear elliptic problems.
\newblock {\em J. Differential Equations}, 275:882--919, 2021.

\bibitem{MR3689156}
D.~Pierotti and G.~Verzini.
\newblock Normalized bound states for the nonlinear {S}chr\"{o}dinger equation
  in bounded domains.
\newblock {\em Calc. Var. Partial Differential Equations}, 56(5):Paper No. 133,
  27, 2017.

\bibitem{pime_voska}
E.~A. Pimentel and V.~Voskanyan.
\newblock Regularity for second-order stationary mean-field games.
\newblock {\em Indiana University Mathematics Journal}, 66(1):1--22, 2017.

\bibitem{Santambrogio2020}
F.~Santambrogio.
\newblock {\em Lecture Notes on Variational Mean Field Games}, pages 159--201.
\newblock Springer International Publishing, Cham, 2020.

\end{thebibliography}
%

\medskip
\small
\begin{flushright}
\noindent \verb"cirant@math.unipd.it"\\
Dipartimento di Matematica ``Tullio Levi-Civita'', Universit\`a di Padova\\ 
Via Trieste 63, 35121 (Italy)
\end{flushright}
\medskip
\small
\begin{flushright}
\noindent \verb"acosenza@math.univ-paris-diderot.fr"\\
Universit\'e Paris Cit\'e and Sorbonne Universit\'e,\\
CNRS, Laboratoire Jacques-Louis Lions (LJLL),\\
F-75006 Paris, France
\end{flushright}
\medskip
\begin{flushright}
\noindent \verb"gianmaria.verzini@polimi.it"\\
Dipartimento di Matematica, Politecnico di Milano\\ 
piazza Leonardo da Vinci 32, 20133 Milano (Italy)
\end{flushright}

\end{document}